\documentstyle[a4,11pt]{article}
\begin{document}
\begin{titlepage}
\vskip 2cm
\begin{flushright}
Preprint CNLP-1996-03
\end{flushright}
\vskip 2cm
\begin{center}
{\bf
Soliton equations in N-dimensions
as exact reductions of Self-dual
Yang - Mills  equation III. Soliton geometry in 2+1 and 3+1 dimensions}
\footnote{Preprint
CNLP-1996-03. Alma-Ata. 1996 }
\end{center}
\vskip 2cm
\begin{center}
Kur. R. Myrzakul and R.  Myrzakulov
\footnote{E-mail: cnlpmyra@satsun.sci.kz}
\end{center}
\vskip 1cm

\begin{center}
 Centre for Nonlinear Problems, PO Box 30, 480035, Alma-Ata-35, Kazakhstan
\end{center}

\begin{abstract}
Some aspects of the relation between differential geometry of curves and
surfaces  and multi-dimensional soliton equations are discussed.
The connection between multi-dimensional soliton equations and
the Self-dual Yang-Mills equation is considered.
\end{abstract}


\end{titlepage}

\setcounter{page}{1}
\newpage

\tableofcontents
\section{Introduction}

One of the multidimensional integrable equations, important both in
physics and mathematics, is the self-dual Yang-Mills (SDYM) equation
[2, 25]. It arises in relatiivity [26, 3] and in field theory [27].
The SDYM equations describe a connection for a bundle over the
Grassmannian of two-dimensional subspaces of the twistor space.
Integrability for a SDYM connection means that its curvature vanishes
on certain two-planes in the tangent space of the Grassmannian. As
shown in [4,5], This allows one to characterize SDYM connections
in terms of the spliting problem for a transition function in a holomorphic
bundle over the Riemann sphere, i.e. the trivialiization of the bundle
[28, 29].

It is well known that the SDYM equation admits  reductions [2, 4-5, 18-19]
(see also the book [1] and references therein) to many
soliton equations in  1+1 and 2+1 dimensions, such as:
i) the Boussinesq, Korteweg-de Vries  (KdV) and modified KdV,
nonlinear Schrodinger, N-wave, Sine-Gordon, Ernst and
chiral field  equations in 1+1 dimensions;
ii) the Kadomtsev - Petviashvili (KP), Davey - Stewartson (DS),
Bogomolny, Toda molecule and the  Ward (chiral field
with torsion) equations,  in 2+1 dimensions.
Moreover, all integrable systems are may be  some reductions of
the SDYM equation or its generalizations (see, e.g. [2, 30]).
It is interesting
to note that the SDYM equation contains also some known ordinary differential
equations (ODE) as one-dimensional reductions [1]. These ODE include: the
Euler-Arnold top, Kovalevskaya  tops, Nahm, Chazy and Painleve  equations.
Also we note that between the Einstein and
Yang-Mills  equations as well exist some connections [3].

On the other hand, it is well known that there exist some another
equations - the mM-LXII and M-LXII equations which contains
many (may be all) soliton equations in 1+1, 2+1 and 3+1 dimensions,
including above mentioned.
These equations are the integrability conditions of the system describing
the moving orthogonal trihedral of a curve or surface [8].

First  aim of the present paper is to show that soliton equations in 2+1 and
3+1 dimensions as well as their "master equations" - the mM-LXII and M-LXII
equations are exact reductions of "premaster equation" - the SDYM equation.
We do this by considering the SDYM equation for GL(Z,C) group (Z=2,3).

The other question which we will discuss below is the relation between
differential geometry (curves, surfaces) and multidimensional soliton
equations. Commonly known fact: between geometry and integrable systems
exists a deeper connection [8-17, 20-22]. Classical examples demonstrating
such connections
are:  the Liouville equation which describes minimals surfaces in $E^{3}$
and the sine-Gordon equation which encodes the whole geometry of
the pseudospherical surfaces, i.e. surfaces with negative constant
Gaussian curvature. The discovery of the Inverse Scattering Transform
(IST) method [31] inspired a lot of modern work on the study of the
connection
between geometry and integrable systems.
So the discussion some aspects of a nature of the tandem "geometry-soliton"
in 2+1 and 3+1 dimensions is second aim of this paper.
We feel our older and permanent love to spin systems (Heisenberg magnet
models = Landau-Lifshitz equations) and below we try consider to some
underlying problems from our "spin" point of view, that is, third aim
of this note.

\section{The (1+1)-dimensional equations of curves/surfaces}

In this section we review briefly the some known basic facts
from 1+1 dimensions to set our notations and terminology.
Consider the curve in 3-dimensional space
or two-dimensional surface.
Equations of such curves and surfaces,
following [8] we can write in the form
$$
\left ( \begin{array}{ccc}
{\bf e}_{1} \\
{\bf e}_{2} \\
{\bf e}_{3}
\end{array} \right)_{x}= C
\left ( \begin{array}{ccc}
{\bf e}_{1} \\
{\bf e}_{2} \\
{\bf e}_{3}
\end{array} \right)
\eqno(1a)
$$
$$
\left ( \begin{array}{ccc}
{\bf e}_{1} \\
{\bf e}_{2} \\
{\bf e}_{3}
\end{array} \right)_{t}= G
\left ( \begin{array}{ccc}
{\bf e}_{1} \\
{\bf e}_{2} \\
{\bf e}_{3}
\end{array} \right) \eqno(1b)
$$
where
$$
C =
\left ( \begin{array}{ccc}
0             & k     & -\sigma \\
-\beta k      & 0     & \tau  \\
\beta \sigma  & -\tau & 0
\end{array} \right) ,\quad
G =
\left ( \begin{array}{ccc}
0       & \omega_{3}  & -\omega_{2} \\
-\beta\omega_{3} & 0      & \omega_{1} \\
\beta\omega_{2}  & -\omega_{1} & 0
\end{array} \right) \eqno(2)
$$
Here ${\bf e}_{j}$ is the moving trihedral of a curve or surface,
$$
{\bf e}_{1}^{2}=\beta = \pm 1, \quad {\bf e}_{2}^{2} = {\bf e}^{2}_{3} = 1
\eqno(3)
$$
$k, \sigma$ and $\tau$ are called the normal curvature, geodesic curvature,
and geodesic torsion, respectevily. So we have
$$
C_{t} - G_{x} + [C, G] = 0 \eqno(4)
$$

In the curves case, the  equation (1a) is as $\sigma = 0$
the Serret - Frenet equation (SFE). In surface interpretation, (1) is the
Gauss-Weingarden equation (GWE) while (4) is the Codazzi-Mainardi-Peterson
equation (CMPE) in orthogonal basis. In the latter case, $k, \tau,
 \sigma, \omega_{j}$ are some functions of $E, F, G, L, M, N $ (coefficients
of first and second fundamental forms).   Let us rewrite the equation (4)
in the form
$$
U_{t}-V_{x}+[U,V]=0  \eqno(5)
$$
where
$$
U =
\frac{1}{2i}\left ( \begin{array}{cc}
\tau        & k -i\sigma \\
k+i\sigma   & -\tau
\end{array} \right) ,\quad
V =
\frac{1}{2i}
\left ( \begin{array}{cc}
\omega_{1}             & \omega_{3}-i\omega_{2} \\
\omega_{3}+i\omega_{2} & - \omega_{1}
\end{array} \right) \eqno(6)
$$
The Lax representation (LR) of the equation (4)-(5) is given by
$$
\psi_{x} = U \psi \eqno(7a)
$$
$$
\psi_{t} =V\psi  \eqno(7b)
$$

Usually, we assume that [8]
$$
k=\sum_{j} \lambda^{j}k_{j}, \quad
\tau=\sum_{j} \lambda^{j}\tau_{j}, \quad
\sigma=\sum_{j} \lambda^{j}\sigma_{j}, \quad
\omega_{k}=\sum_{j} \lambda^{j}\omega_{kj} \eqno(8)
$$
where $j=0, \pm 1, \pm 2, \pm 3, ...$ .  Sometimes, instead of (8) we take
the following more general case
$$
k=\sum_{j} h_{1j}k_{j}, \quad
\tau=\sum_{j} h_{2j}\tau_{j}, \quad
\sigma=\sum_{j} h_{3j}\sigma_{j}, \quad
\omega_{k}=\sum_{j} h_{4j}\omega_{kj} \eqno(9)
$$
Here $h_{ij}= h_{ij}(\lambda)$ and   $\lambda$ is some characteristic
parameter of curves or surfaces or some  function of such parameters.
In the cases (8) and (9), i.e when
$k, \tau, \sigma, \omega_{j}$ are some functions of $\lambda$,
the equations (1), (4)=(5) and  (7),  we call the M-LXIX, M-LXX
and M-LXXI equations respectively.
[ These  conditional notations we use in order to accurately distinguish
these equations from the case when $k, \tau, \sigma, \omega_{j}$ are
independent of $\lambda$ and for convenience in our internal working
kitchen. The more so, as of course it is not means that
we have e.g. 72 equations.]  Particular cases.

i) {\it The Gauss-Weingarden  and Codazzi-Mainardi-Peterson equations}.
These equations  correspond to the case
$$
k=k_{0}, \quad
\tau=\tau_{0}, \quad
\sigma=\sigma_{0}, \quad
\omega_{k}=\omega_{k0} \eqno(10)
$$

ii) {\it The ZS-AKNS problem}.
The famous Zakharov-
Shabat-Ablovitz-Kaup-Newell-Segur (ZS-AKNS)  spectral
problem which generate many soliton equations in 1+1  is the
particular case of the M-LXXI equation (7a) as
$$
k_{0}=i(p+q), \quad k_{j}=0, \quad j\neq 0 \eqno(11a)
$$
$$
\sigma_{0} = p-q, \quad \sigma_{j} = 0, \quad j \neq 0  \eqno(11b)
$$
$$
\tau_{1} = -2, \quad \tau_{j} = 0, \quad j \neq 1  \eqno(11c)
$$

iii) {\it The Kaup-Newell-Wadati-Konno-Ichicawa spectral problem}.
This case corresponds to reduction
$$
k_{1}=i\lambda (p+q), \quad k_{j}=0, \quad j \neq 1 \eqno(11d)
$$
$$
\sigma_{1} = \lambda (p-q), \quad \sigma_{j} = 0, \quad j \neq 1  \eqno(11e)
$$
$$
\tau_{1} = -2, \quad \tau_{j} = 0, \quad j \neq 1  \eqno(11f)
$$

iv) {\it The  equations of a principal chiral field}.
The equations of a principal chiral field for functions $u, v$
$$
u_{t}+\frac{1}{2}[u,v]=0, \quad v_{x}-\frac{1}{2}[u,v]=0 \eqno(11g)
$$
can equally be represented with the aid of (4) or (5), if we choose
in this case
$$
h_{l1}=\frac{1}{1-\lambda}, \quad k_{j}=\tau_{j}=\sigma_{j}=0, \quad j \neq 1
\eqno(11h)
$$
$$
U=\frac{u}{1-\lambda}, \quad V=\frac{v}{\lambda + 1} \eqno(11i)
$$

v) {\it The (1+1)-dimensional mM-LXVI equation}.
This case corresponds to the reduction when
$$
k^{2}+\tau^{2}+\sigma^{2} = n^{2} \eqno(11j)
$$

So we see that practically all (1+1)-dimensional soliton equations
can be obtained from the M-LXX equation (4)=(5) as some reductions.

\section{Curves and Surfaces in 4-dimensions: the mM-LXI equation}

\subsection{The (3+1)-dimensional mM-LXI equation: the space case}

The (3+1)-dimensional mM-LXI equation reads as [8]
$$
\left ( \begin{array}{c}
{\bf e}_{1} \\
{\bf e}_{2} \\
{\bf e}_{3}
\end{array} \right)_{x}= A
\left ( \begin{array}{ccc}
{\bf e}_{1} \\
{\bf e}_{2} \\
{\bf e}_{3}
\end{array} \right), \quad \left ( \begin{array}{ccc}
{\bf e}_{1} \\
{\bf e}_{2} \\
{\bf e}_{3}
\end{array} \right)_{y}= B
\left ( \begin{array}{ccc}
{\bf e}_{1} \\
{\bf e}_{2} \\
{\bf e}_{3}
\end{array} \right),
$$
$$
\left ( \begin{array}{ccc}
{\bf e}_{1} \\
{\bf e}_{2} \\
{\bf e}_{3}
\end{array} \right)_{z}= C
\left ( \begin{array}{ccc}
{\bf e}_{1} \\
{\bf e}_{2} \\
{\bf e}_{3}
\end{array} \right), \quad \left ( \begin{array}{ccc}
{\bf e}_{1} \\
{\bf e}_{2} \\
{\bf e}_{3}
\end{array} \right)_{t}= D
\left ( \begin{array}{ccc}
{\bf e}_{1} \\
{\bf e}_{2} \\
{\bf e}_{3}
\end{array} \right)
\eqno(12)
$$
where
$$
A =
\left ( \begin{array}{ccc}
0             & k     &-\sigma \\
-\beta k      & 0     & \tau  \\
\beta\sigma   & -\tau & 0
\end{array} \right) ,
\quad
B=
\left ( \begin{array}{ccc}
0            & m_{3}  & -m_{2} \\
-\beta m_{3} & 0      & m_{1} \\
\beta m_{2}  & -m_{1} & 0
\end{array} \right)
$$
$$
C=
\left ( \begin{array}{ccc}
0            & n_{3}  & -n_{2} \\
-\beta n_{3} & 0      & n_{1} \\
\beta n_{2}  & -n_{1} & 0
\end{array} \right),
\quad
D=
\left ( \begin{array}{ccc}
0            & \omega_{3}  & -\omega_{2} \\
-\beta \omega_{3} & 0      & \omega_{1} \\
\beta \omega_{2}  & -\omega_{1} & 0
\end{array} \right).  \eqno(13)
$$

\subsection{The (3+1)-dimensional modified M-LXII equation: the space case}

From (12), we obtain the following mM-LXII equation [8]
$$
A_{y} - B_{x} + [A,B] = 0 \eqno (14a)
$$
$$
A_{z} - C_{x} + [A,C] = 0 \eqno (14b)
$$
$$
A_{t} - D_{x} + [A,D] = 0 \eqno (14c)
$$
$$
B_{z} - C_{y} + [B,C] = 0 \eqno (14d)
$$
$$
B_{t} - D_{y} + [B,D] = 0 \eqno (14e)
$$
$$
C_{t} - D_{z} + [C,D] = 0 \eqno (14f)
$$

\subsubsection{Lax representation of the (3+1)-dimensional mM-LXII equation}

The mM-LXII equation (14) we can rewrite in the following form
$$
F_{xy} = U_{y}-V_{x} + [U,V]=0  \eqno(15a)
$$
$$
F_{xz} = U_{z}-W_{x} + [U,W]=0  \eqno(15b)
$$
$$
F_{xt}=U_{t}-T_{x} + [U,T]=0  \eqno(15c)
$$
$$
F_{yz}=V_{z}-W_{y} + [V,W]=0  \eqno(15d)
$$
$$
F_{yt}=V_{t}-T_{y} + [V,T]=0  \eqno(15e)
$$
$$
F_{zt}=W_{t}-T_{z} + [W,T]=0  \eqno(15f)
$$
where
$$
U =
\frac{1}{2i}\left ( \begin{array}{cc}
\tau        & k -i\sigma \\
k+i\sigma   & -\tau
\end{array} \right) ,\quad
T =
\frac{1}{2i}
\left ( \begin{array}{cc}
\omega_{1}             & \omega_{3}-i\omega_{2} \\
\omega_{3}+i\omega_{2} & - \omega_{1}
\end{array} \right)
$$
$$
V =
\frac{1}{2i}\left ( \begin{array}{cc}
m_{1}        & m_{3} -im_{2} \\
m_{3}+im_{2}   & -m_{1}
\end{array} \right) ,\quad
W =
\frac{1}{2i}
\left ( \begin{array}{cc}
n_{1}             & n_{3}-in_{2} \\
n_{3}+in_{2} & - n_{1}
\end{array} \right) \eqno(16)
$$
The LR of the mM-LXII equation (14) or (15) is given by
$$
g_{x} = Ug, \quad g_{y} = Vg, \quad g_{z} = Wg, \quad g_{t} = Tg \eqno(17)
$$
We note that this LR can be rewritten in the more usual form
$$
Lg = 0, \quad Mg = 0 \eqno(18)
$$
where, for example, $L, M$ we can take in the form
$$
L = \partial_{x}+a\lambda \partial_{y} +b\lambda^{2} \partial_{z} -
(U+a\lambda V +b W), \eqno(19a)
$$
$$
M=\partial_{t}+e\lambda^{3}\partial_{y}+f\lambda^{4}\partial_{z}
-(T+e\lambda^{3}V+f\lambda^{4}W) \eqno(19b)
$$
or
$$
L = D_{x}+a\lambda D_{y}+b\lambda^{2} D_{z}, \quad M=D_{t}+e\lambda^{3}D_{y}
+f\lambda^{4}D_{z} \eqno(20)
$$
Here $D_{j}$ are the covariant derivatives
$$
D_{x} = \partial_{x}-U, \quad
D_{y} = \partial_{y}-V, \quad
D_{z} = \partial_{z}-W, \quad
D_{t} = \partial_{t}-T  \eqno(21)
$$

\subsection{The (3+1)-dimensional mM-LXI equation: the plane case}

The mM-LXI equation in plane has the form [8]
$$
\left ( \begin{array}{cc}
{\bf e}_{1} \\
{\bf e}_{2}
\end{array} \right)_{x}= A_{p}
\left ( \begin{array}{cc}
{\bf e}_{1} \\
{\bf e}_{2}
\end{array} \right), \quad
\left ( \begin{array}{cc}
{\bf e}_{1} \\
{\bf e}_{2}
\end{array} \right)_{y}= B_{p}
\left ( \begin{array}{cc}
{\bf e}_{1} \\
{\bf e}_{2}
\end{array} \right),
$$
$$
\left ( \begin{array}{cc}
{\bf e}_{1} \\
{\bf e}_{2}
\end{array} \right)_{z}= C_{p}
\left ( \begin{array}{cc}
{\bf e}_{1} \\
{\bf e}_{2}
\end{array} \right), \quad
\left ( \begin{array}{cc}
{\bf e}_{1} \\
{\bf e}_{2}
\end{array} \right)_{t}= D_{p}
\left ( \begin{array}{cc}
{\bf e}_{1} \\
{\bf e}_{2}
\end{array} \right)
\eqno(22)
$$
where
$$
A_{p} =
\left ( \begin{array}{cc}
0             & k     \\
-\beta k      & 0
\end{array} \right) ,
\quad
B_{p}=
\left ( \begin{array}{cc}
0            & m_{3}   \\
-\beta m_{3} & 0
\end{array} \right)
$$
$$
C_{p} =
\left ( \begin{array}{cc}
0             & n_{3}     \\
-\beta n_{3}      & 0
\end{array} \right) ,
\quad
D_{p}=
\left ( \begin{array}{cc}
0            & \omega_{3}   \\
-\beta \omega_{3} & 0
\end{array} \right).  \eqno(23)
$$

\subsection{The (3+1)-dimensional modified M-LXII equation: the plane case}

In the plane case the mM-LXII equation takes the following simple form
$$
k_{y} = m_{3x}, \quad k_{z} = n_{3x}, \quad m_{3z} = n_{3y},
\quad m_{3t}=\omega_{3y}, \quad n_{3t} = \omega_{3z} \eqno(24a)
$$
Hence we get
$$
m_{3} = \partial^{-1}_{x}k_{y}, \quad n_{3} =\partial^{-1}_{x}k_{z}
\eqno (24b)
$$
The (3+1)-dimensional nonlinear evolution equation for $k$ has the form
$$
k_{t}= \omega_{3x}  \eqno(25)
$$
We note that this equation gives the (3+1)-dimensional KdV, KP, NV, mNV and
so on,  equations as some reductions with the corresponding LR (17).

\subsection{The  M-LXI and M-LXII equations}

The  M-LXI and M-LXII equations are the particular cases of the mM-LXI
and mM-LXII equations as $\sigma = 0$ respectively.

\section{Curves and Surfaces in 4 dimensions: the modified M-LXVIII equation}

\subsection{The modified M-LXVIII equation: the space case }

One of important four dimensional  equation of curves and/or surfaces
is the following modified M-LXVIII (mM-LXVIII) equation
$$
\left ( \begin{array}{ccc}
{\bf e}_{1} \\
{\bf e}_{2} \\
{\bf e}_{3}
\end{array} \right)_{\xi_{1}}=a
\left ( \begin{array}{ccc}
{\bf e}_{1} \\
{\bf e}_{2} \\
{\bf e}_{3}
\end{array} \right)_{\xi_{3}}+C
\left ( \begin{array}{ccc}
{\bf e}_{1} \\
{\bf e}_{2} \\
{\bf e}_{3}
\end{array} \right)
\eqno(26a)
$$
$$
\left ( \begin{array}{ccc}
{\bf e}_{1} \\
{\bf e}_{2} \\
{\bf e}_{3}
\end{array} \right)_{\xi_{2}}=b
\left ( \begin{array}{ccc}
{\bf e}_{1} \\
{\bf e}_{2} \\
{\bf e}_{3}
\end{array} \right)_{\xi_{4}}+ G
\left ( \begin{array}{ccc}
{\bf e}_{1} \\
{\bf e}_{2} \\
{\bf e}_{3}
\end{array} \right) \eqno(26b)
$$
where $\xi_{j}$ are may be real or complex coordinates, $a, b$
are some parameters,
$$
C =
\left ( \begin{array}{ccc}
0             & k     & -\sigma \\
-\beta k      & 0     & \tau  \\
\beta \sigma  & -\tau & 0
\end{array} \right) ,\quad
G =
\left ( \begin{array}{ccc}
0       & \omega_{3}  & -\omega_{2} \\
-\beta\omega_{3} & 0      & \omega_{1} \\
\beta\omega_{2}  & -\omega_{1} & 0
\end{array} \right) \eqno(27)
$$
The compativility condition of the equations (26) is given by
$$
C_{0\xi_{2}}-bC_{0\xi_{4}} +aG_{1\xi_{3}}-G_{1\xi_{1}}+[C_{0},G_{1}]=0
\eqno(28)
$$

The LR of this equation can be written as
$$
\psi_{\xi_{1}}= a\psi_{\xi_{3}} + B_{0}g \eqno(29a)
$$
$$
\psi_{\xi_{2}}= b\psi_{\xi_{4}} + B_{1}g \eqno(29b)
$$
Hence we get
$$
B_{0\xi_{2}}-bB_{0\xi_{4}} +aB_{1\xi_{3}}-B_{1\xi_{1}}+[B_{0},B_{1}]=0
\eqno(30)
$$
with
$$
B_{0} =
\frac{1}{2i}\left ( \begin{array}{cc}
\tau        & k -i\sigma \\
k+i\sigma   & -\tau
\end{array} \right) ,\quad
B_{1}=
\frac{1}{2i}
\left ( \begin{array}{cc}
\omega_{1}             & \omega_{3}-i\omega_{2} \\
\omega_{3}+i\omega_{2} & - \omega_{1}
\end{array} \right) \eqno(31)
$$

\subsection{The mM-LXVIII equation: the plane case }

The mM-LXVIII equation in the plane case  has the form
$$
\left ( \begin{array}{cc}
{\bf e}_{1} \\
{\bf e}_{2}
\end{array} \right)_{\xi_{1}}=a
\left ( \begin{array}{cc}
{\bf e}_{1} \\
{\bf e}_{2}
\end{array} \right)_{\xi_{3}}+C
\left ( \begin{array}{cc}
{\bf e}_{1} \\
{\bf e}_{2}
\end{array} \right)
\eqno(32a)
$$
$$
\left ( \begin{array}{cc}
{\bf e}_{1} \\
{\bf e}_{2}
\end{array} \right)_{\xi_{2}}=b
\left ( \begin{array}{cc}
{\bf e}_{1} \\
{\bf e}_{2}
\end{array} \right)_{\xi_{4}}+ G
\left ( \begin{array}{cc}
{\bf e}_{1} \\
{\bf e}_{2}
\end{array} \right) \eqno(32b)
$$

\subsection{The M-LXVIII equation}

The M-LXVIII equation is the particular case of the mM-LXVIII equation
as $\sigma = 0$.

\section{SDYM equation}

It is standard to define
a Yang-Mills vector  bundle over a four-dimensional manifold $M$
with connection one-form $A=A_{\mu}(x^{\nu})dx^{\mu}$.
The SDYM equations in this manifolds let us  write to set our notation:
$$
F = \ast F  \eqno(33)
$$
where $F$ is a curvature 2-form pulled back to M from the gauge bundle
$P(M, g)$, explicitly:
$$
F = d\omega + \omega \wedge \omega  \eqno(34)
$$
Here the connection 1-form $\omega$ on $P$ takes values in the Lie
algebra $g$ of the gauge group $G$. In terms of Cartesian coordinates
$x^{\mu}$, they can be expressed as
$$
F_{\mu\nu} = \frac{1}{2}\epsilon_{\mu\nu\kappa\sigma}F_{\kappa\sigma}
\eqno(35)
$$
where $\mu, \nu, ... =1,2,3,4,  \quad \epsilon_{\mu\nu\kappa\sigma}$ stands
for the completely antisymmetric
tensor in four dimensions with the convention: $\epsilon_{1234} = 1$.  The
components of the field strength $(F_{\mu\nu})$ are given by
$$
F_{\mu\nu} = \partial_{\mu}A_{\nu} - \partial_{\nu}A_{\mu}
+ [A_{\mu}, A_{\nu}]\eqno(36)
$$
We note that the SDYM equations (33) are invariant under the gauge transformation
$$
A_{\mu} \rightarrow \phi^{-1} A_{\mu} \phi - \phi^{-1}\phi_{\mu} \eqno(37)
$$

Let us introduce the some null coordinates $\xi_{j}$ in which the metric on $M$
is  $ds^{2}=d\xi_{1} d\xi_{3} + d\xi_{2} d\xi_{4}$.  For these coordinates,
the SDYM equations are given by
$$
F_{\xi_{1}\xi_{2}} = 0  \eqno(38a)
$$
$$
F_{\xi_{3}\xi_{4}} = 0  \eqno(38b)
$$
$$
F_{\xi_{1}\xi_{4}} - F_{\xi_{2}\xi_{3}} = 0  \eqno(38c)
$$
or
$$
\partial_{\xi_{1}}A_{2}-\partial_{\xi_{2}}A_{1}
+[A_{1}, A_{2}] = 0 \eqno(39a)
$$
$$
\partial_{\xi_{3}}A_{4}-\partial_{\xi_{4}}A_{3}
+[A_{3}, A_{4}] = 0 \eqno(39b)
$$
$$
\partial_{\xi_{1}}A_{4}-\partial_{\xi_{4}}A_{1}
+[A_{1}, A_{4}] =
\partial_{\xi_{2}}A_{3}-\partial_{\xi_{3}}A_{2}
+[A_{2}, A_{3}] \eqno(39c)
$$
For the  equations (38)-(39), the LR has the form
$$
L\Phi = 0,  \quad M\Phi = 0  \eqno(40)
$$
where
$$
L = D_{\xi_{1}} - \lambda D_{\xi_{3}}, \quad M= D_{\xi_{2}}-
\lambda D_{\xi_{4}}  \eqno(41)
$$
and $\lambda$ is the spectral parameter.

\section{The 4-dimensional mM-LXVIII equation and SDYM equation}

To establish the connection between the 4-dimensional mM-LXVIII
equation and SDYM equation,  we rewrite the LR (29) of the mM-LXVIII
equation in the form
$$
L\psi =0, \quad  M\psi = 0  \eqno(42)
$$
where
$$
L = D_{\xi_{1}} - a D_{\xi_{3}}, \quad   M = D_{\xi_{2}} - b D_{\xi_{4}}
\eqno(43)
$$
Here
$$
D_{\xi_{j}} = \partial_{\xi_{j}} - A_{j}, \quad A_{1}-aA_{3}=B_{0},
\quad A_{2}-bA_{4} = B_{1}   \eqno(44)
$$
We see that the equations (40) and (42) have the same form if $a=b=\lambda$.
So the 4-dimensional mM-LXVIII equation (28)=(30) is equivalent to the
 SDYM equations (38)=(39) as $a=b=\lambda$. Hence follows that the curve or
surfaces (26) is integrable in this case and in this sense, i.e. they
are the integrable curve or integrable surface in four (4=3+1=2+2=1+3)
dimensions.

\section{The 4-dimensional  mM-LXII equations and SDYM equation}

To establish the connection between the mM-LXII equation (14)=(15) and
the SDYM equation
(38)=(39) we consider the following transformation for coordinates
$$
\left ( \begin{array}{cccc}
\xi_{1} \\
\xi_{2} \\
\xi_{3} \\
\xi_{4}
\end{array} \right) = H
\left ( \begin{array}{cccc}
x \\
y \\
z \\
t
\end{array} \right) \eqno(45)
$$
where
$$
H=
\left ( \begin{array}{cccc}
a_{11}  & a_{12} & a_{13} & a_{14} \\
a_{21}  & a_{22} & a_{23} & a_{24} \\
a_{31}  & a_{32} & a_{33} & a_{34} \\
a_{41}  & a_{42} & a_{43} & a_{44}
\end{array} \right) \eqno(46)
$$
Also we use the expressions for $U,V,W,T$
$$
U=b_{11}A_{\xi_{1}} + b_{12}A_{\xi_{2}} + b_{13}A_{\xi_{3}}
+ b_{14}A_{\xi_{4}} \eqno(47a)
$$
$$
V=b_{21}A_{\xi_{1}} + b_{22}A_{\xi_{2}} + b_{23}A_{\xi_{3}}
+ b_{24}A_{\xi_{4}} \eqno(47b)
$$
$$
W=b_{31}A_{\xi_{1}} + b_{32}A_{\xi_{2}} + b_{33}A_{\xi_{3}}
+ b_{34}A_{\xi_{4}} \eqno(47c)
$$
$$
T=b_{41}A_{\xi_{1}} + b_{42}A_{\xi_{2}} + b_{43}A_{\xi_{3}}
+ b_{44}A_{\xi_{4}} \eqno(47d)
$$
We can consider the other may be more general transformation
$$
\xi_{j} = f_{j}(x,y,z,t) \eqno(48)
$$
where $f_{j}$ are some functions.
Then after some algebra we have
$$
F_{xy} = l_{12}^{1}F_{\xi_{1}\xi_{3}} + l_{12}^{2}F_{\xi_{2}\xi_{4}}
+ l_{12}^{3}(F_{\xi_{1}\xi_{2}} -  F_{\xi_{3}\xi_{4}}) \eqno(49a)
$$
$$
F_{xz} = l_{13}^{1}F_{\xi_{1}\xi_{3}} + l_{13}^{2}F_{\xi_{2}\xi_{4}}
+ l_{13}^{3}(F_{\xi_{1}\xi_{2}} -  F_{\xi_{3}\xi_{4}}) \eqno(49b)
$$
$$
F_{xt} = l_{14}^{1}F_{\xi_{1}\xi_{3}} + l_{14}^{2}F_{\xi_{2}\xi_{4}}
+ l_{14}^{3}(F_{\xi_{1}\xi_{2}} -  F_{\xi_{3}\xi_{4}}) \eqno(49c)
$$
$$
F_{yz} = l_{23}^{1}F_{\xi_{1}\xi_{3}} + l_{23}^{2}F_{\xi_{2}\xi_{4}}
+ l_{23}^{3}(F_{\xi_{1}\xi_{2}} -  F_{\xi_{3}\xi_{4}}) \eqno(49d)
$$
$$
F_{yt} = l_{24}^{1}F_{\xi_{1}\xi_{3}} + l_{24}^{2}F_{\xi_{2}\xi_{4}}
+ l_{24}^{3}(F_{\xi_{1}\xi_{2}} -  F_{\xi_{3}\xi_{4}}) \eqno(49e)
$$
$$
F_{zt} = l_{34}^{1}F_{\xi_{1}\xi_{3}} + l_{34}^{2}F_{\xi_{2}\xi_{4}}
+ l_{34}^{3}(F_{\xi_{1}\xi_{2}} -  F_{\xi_{3}\xi_{4}}) \eqno(49f)
$$

So if $F_{\mu\nu}  (\mu,\nu=\xi_{1}, \xi_{2}, \xi_{3}, \xi_{4})$ satisfy
the SDYM equation (38) then $F_{ij} (i,j=x,y,z,t)$ satisfy the mM-LXII
equation (15).

\section{The (2+1)-dimensional mM-LXVIII equation and SDYM equation}

The mM-LXVIII equation in 3=2+1 dimensions can be written for example
in the form
$$
\left ( \begin{array}{ccc}
{\bf e}_{1} \\
{\bf e}_{2} \\
{\bf e}_{3}
\end{array} \right)_{\xi_{1}}=C
\left ( \begin{array}{ccc}
{\bf e}_{1} \\
{\bf e}_{2} \\
{\bf e}_{3}
\end{array} \right)
\eqno(50a)
$$
$$
\left ( \begin{array}{ccc}
{\bf e}_{1} \\
{\bf e}_{2} \\
{\bf e}_{3}
\end{array} \right)_{\xi_{2}}=b
\left ( \begin{array}{ccc}
{\bf e}_{1} \\
{\bf e}_{2} \\
{\bf e}_{3}
\end{array} \right)_{\xi_{4}}+ G
\left ( \begin{array}{ccc}
{\bf e}_{1} \\
{\bf e}_{2} \\
{\bf e}_{3}
\end{array} \right) \eqno(50b)
$$
In this case the corresponding LR takes the form
$$
\psi_{\xi_{1}}=  B_{0}g \eqno(51a)
$$
$$
\psi_{\xi_{2}}= b\psi_{\xi_{4}} + B_{1}g \eqno(51b)
$$
Hence we get
$$
B_{0\xi_{2}}-bB_{0\xi_{4}} -B_{1\xi_{1}}+[B_{0},B_{1}]=0
\eqno(52)
$$

The corresponding (2+1)-dimensional SDYM equations have the forms
$$
\partial_{\xi_{1}}A_{2}-\partial_{\xi_{2}}A_{1}
+[A_{1}, A_{2}] = 0 \eqno(53a)
$$
$$
\partial_{\xi_{4}}A_{3} - [A_{3}, A_{4}] = 0 \eqno(53b)
$$
$$
\partial_{\xi_{1}}A_{4}-\partial_{\xi_{4}}A_{1}
+[A_{1}, A_{4}] =
\partial_{\xi_{2}}A_{3}+[A_{2}, A_{3}] \eqno(53c)
$$

\section{The (2+1)-dimensional mM-LXI and mM-LXII equations: the space case
and  SDYM equation}

\subsection{The (2+1)-dimensional mM-LXI equation for the space case}

The (2+1)-dimensional mM-LXI equation reads as
$$
\left ( \begin{array}{c}
{\bf e}_{1} \\
{\bf e}_{2} \\
{\bf e}_{3}
\end{array} \right)_{x}= A
\left ( \begin{array}{ccc}
{\bf e}_{1} \\
{\bf e}_{2} \\
{\bf e}_{3}
\end{array} \right) \eqno(54a)
$$
$$
\left ( \begin{array}{ccc}
{\bf e}_{1} \\
{\bf e}_{2} \\
{\bf e}_{3}
\end{array} \right)_{y}= B
\left ( \begin{array}{ccc}
{\bf e}_{1} \\
{\bf e}_{2} \\
{\bf e}_{3}
\end{array} \right) \eqno(54b)
$$
$$
\left ( \begin{array}{ccc}
{\bf e}_{1} \\
{\bf e}_{2} \\
{\bf e}_{3}
\end{array} \right)_{t}= D
\left ( \begin{array}{ccc}
{\bf e}_{1} \\
{\bf e}_{2} \\
{\bf e}_{3}
\end{array} \right)
\eqno(54c)
$$

\subsection{The (2+1)-dimensional mM-LXII equation: the space case}

From (54), we obtain the following mM-LXII equation [8]
$$
A_{y} - B_{x} + [A,B] = 0 \eqno (55a)
$$
$$
A_{t} - D_{x} + [A,D] = 0 \eqno (55b)
$$
$$
B_{t} - D_{y} + [B,D] = 0 \eqno (55c)
$$
Some (2+1)-dimensional soliton equations such as:
DS, Zakharov, (2+1)-complex mKdV, (2+1)-dNLS and so on, are exact
reductions of the mM-LXII equation (55).

\subsection{Lax representation of the (2+1)-dimensional mM-LXII equation}

As above the mM-LXII equation (55) we write in the following form
$$
F_{xy} = U_{y}-V_{x} + [U,V]=0  \eqno(56a)
$$
$$
F_{xt}=U_{t}-T_{x} + [U,T]=0  \eqno(56b)
$$
$$
F_{yt}=V_{t}-T_{y} + [V,T]=0  \eqno(56c)
$$

So that the LR of the mM-LXVIII equation (55)=(56) is given by
$$
g_{x} = Ug, \quad g_{y} = Vg,  \quad g_{t} = Tg \eqno(57)
$$
We note that this LR we can write  in the more usual form
$$
Lg = 0, \quad Mg = 0 \eqno(58)
$$
where, for example, $L, M$ we can take in the form
$$
L = \partial_{x}+a\lambda \partial_{y}  -
(U+a\lambda V), \eqno(59a)
$$
$$
M=\partial_{t}+e\lambda^{2}\partial_{y}
-(T+e\lambda^{2}V) \eqno(59b)
$$
or
$$
L = D_{x}+a\lambda D_{y}, \quad M=D_{t}+e\lambda^{3}D_{y}. \eqno(60)
$$

\subsection{The (2+1)-dimensional mM-LXII equation and SDYM equation}

To establish the connection between the mM-LXII equation (55)=(56) and
the SDYM equation
(38) we consider as above the following transformation for coordinates
$$
\left ( \begin{array}{cccc}
\xi_{1} \\
\xi_{2} \\
\xi_{3} \\
\xi_{4}
\end{array} \right) = H
\left ( \begin{array}{cccc}
x \\
y \\
0 \\
t
\end{array} \right) \eqno(61)
$$
where
$$
H=
\left ( \begin{array}{cccc}
a_{11}  & a_{12} & a_{13} & a_{14} \\
a_{21}  & a_{22} & a_{23} & a_{24} \\
a_{31}  & a_{32} & a_{33} & a_{34} \\
a_{41}  & a_{42} & a_{43} & a_{44}
\end{array} \right) \eqno(62)
$$
Also we use the expressions for $U,V,W,T$
$$
U=b_{11}A_{\xi_{1}} + b_{12}A_{\xi_{2}} + b_{13}A_{\xi_{3}}
+ b_{14}A_{\xi_{4}} \eqno(63a)
$$
$$
V=b_{21}A_{\xi_{1}} + b_{22}A_{\xi_{2}} + b_{23}A_{\xi_{3}}
+ b_{24}A_{\xi_{4}} \eqno(63b)
$$
$$
T=b_{41}A_{\xi_{1}} + b_{42}A_{\xi_{2}} + b_{43}A_{\xi_{3}}
+ b_{44}A_{\xi_{4}} \eqno(63c)
$$
We can consider the other may be more general transformation for coordinates
$$
\xi_{j} = f_{j}(x,y,t) \eqno(64)
$$

Then after some algebra we have
$$
F_{xy} = l_{12}^{1}F_{\xi_{1}\xi_{3}} + l_{12}^{2}F_{\xi_{2}\xi_{4}}
+ l_{12}^{3}(F_{\xi_{1}\xi_{2}} -  F_{\xi_{3}\xi_{4}}) \eqno(65a)
$$
$$
F_{xt} = l_{14}^{1}F_{\xi_{1}\xi_{3}} + l_{14}^{2}F_{\xi_{2}\xi_{4}}
+ l_{14}^{3}(F_{\xi_{1}\xi_{2}} -  F_{\xi_{3}\xi_{4}}) \eqno(65b)
$$
$$
F_{yt} = l_{24}^{1}F_{\xi_{1}\xi_{3}} + l_{24}^{2}F_{\xi_{2}\xi_{4}}
+ l_{24}^{3}(F_{\xi_{1}\xi_{2}} -  F_{\xi_{3}\xi_{4}}) \eqno(65c)
$$

So if $F_{\mu\nu}  (\mu,\nu=\xi_{1}, \xi_{2}, \xi_{3}, \xi_{4})$ satisfy
the SDYM equation (38) then $F_{ij} (i,j=x,y,t)$ satisfy the mM-LXII
equation (55)=(56).

\section{The (2+1)-dimensional mM-LXI and mM-LXII equations: the plane case}

\subsection{The (2+1)-dimensional mM-LXI equation: the plane case}

The (2+1)-dimensional mM-LXI equation in plane has the form [8]
$$
\left ( \begin{array}{cc}
{\bf e}_{1} \\
{\bf e}_{2}
\end{array} \right)_{x}= A_{p}
\left ( \begin{array}{cc}
{\bf e}_{1} \\
{\bf e}_{2}
\end{array} \right), \quad
\left ( \begin{array}{cc}
{\bf e}_{1} \\
{\bf e}_{2}
\end{array} \right)_{y}= B_{p}
\left ( \begin{array}{cc}
{\bf e}_{1} \\
{\bf e}_{2}
\end{array} \right), \quad
\left ( \begin{array}{cc}
{\bf e}_{1} \\
{\bf e}_{2}
\end{array} \right)_{t}= D_{p}
\left ( \begin{array}{cc}
{\bf e}_{1} \\
{\bf e}_{2}
\end{array} \right)
\eqno(66)
$$

\subsection{The (2+1)-dimensional m M-LXII equation: the plane case}

In the plane case the mM-LXII equation takes the following simple form
$$
k_{y} = m_{3x}, \quad m_{3t}=\omega_{y}  \eqno(67a)
$$
Hence we get
$$
m_{3} = \partial^{-1}_{y}k_{y}\eqno (67b)
$$
The nonlinear evolution equation has the form
$$
k_{t}= \omega_{3x}  \eqno(68)
$$
Many (2+1)-dimensional integrable equations such as the
Kadomtsev-Petviashvili, Novikov-Veselov (NV), mNV, KNV, (2+1)-KdV, mKdV
equations are the integrable reductions of
the M-LXII equation (68).

\section{The  (2+1)-dimensional M-LXI and M-LXII equations}

The  (2+1)-dimensional M-LXI and M-LXII equations are the particular cases
of the (2+1)-dimensional mM-LXI
and mM-LXII equations as $\sigma = 0$ respectively.

\section{The (2+1)-dimensional mM-LXII equation as the reduction of SDYM
equation}

One of a way to establish the connection between the mM-LXII
equation (55)-(56), i.e. soliton equations in 2+1  and the SDYM
equation (38) is the following.
Consider the coordinates
$$
\xi_{1} = it, \quad \xi_{2} = -it, \quad \xi_{3} = x+iy, \quad
\xi_{4} = x-iy, \eqno(69)
$$
Now in the SDYM equation (38)=(39) we  take
$$
A_{1}=-iD, \quad A_{2}=iD, \quad A_{3}=A-iB, \quad A_{4}=A+iB.  \eqno(70)
$$
where we mention that $A,B,C,D$  are in   our case real matrices.
Then the SDYM equation (39) reduces to the (2+1)-dimensional
mM-LXII equation (56).
So the (2+1)-dimensional mM-LXII equation is the integrable reduction
of the SDYM equation.

\section{Soliton equations in 2+1 dimensions are exact reductions
of  SDYM equation}

As many (may be all) soliton equations in 2+1 dimensions
are some integrable reductions of the mM-LXII and/or M-LXII equations
(55) and/or (56)  then as follows
from the results of the previous section these (2+1)-dimensional
soliton equations are exact reductions of the SDYM equation (see, e.g
[18, 22]).

\section{The mM-LXII and Soliton equations in 3+1 dimensions as
exact reductions of  SDYM equation}

Here we have the same arguments
as in sections 12 and 13 [19, 22].

\section{Surfaces with the restriction: $k^{2}+\tau^{2}+
\sigma^{2} = n^{2}(x,y,z,t)$.
The M-LXVI  equation  and spin systems. L-equivalence and G-equivalence}

Until now we not spoken about spin systems although  not forget their.
They are to us something greater than only mathematical object.
So that let us find the place for them in our formalism and in this note.
To this purpose, we remember the M-LXVI equation which is the particular
case of some above considered equations as
$$
k^{2} + \tau^{2} + \sigma^{2} =  n^{2}(x,y,z,t)  \eqno(71)
$$
Below we consider the case when $n=constant$. Let
$$
k = n S_{1}, \quad \sigma = n S_{2}, \tau = n S_{3}  \eqno(72)
$$
Then from (71) follows that
$$
S_{1}^{2} + S_{2}^{2} + S_{3}^{2} = 1 \eqno(73)
$$
Consider the simplest example,  the (1+1)-dimensional M-LXVI equation
having the form
$$
\left ( \begin{array}{ccc}
{\bf f}_{1} \\
{\bf f}_{2} \\
{\bf f}_{3}
\end{array} \right)_{x}= C^{\prime}
\left ( \begin{array}{ccc}
{\bf f}_{1} \\
{\bf f}_{2} \\
{\bf f}_{3}
\end{array} \right)
\eqno(74a)
$$
$$
\left ( \begin{array}{ccc}
{\bf f}_{1} \\
{\bf f}_{2} \\
{\bf f}_{3}
\end{array} \right)_{t}= G^{\prime}
\left ( \begin{array}{ccc}
{\bf f}_{1} \\
{\bf f}_{2} \\
{\bf f}_{3}
\end{array} \right) \eqno(74b)
$$
where
$$
C^{\prime} =n
\left ( \begin{array}{ccc}
0             & S_{1}     & -S_{2} \\
-\beta S_{1}      & 0     & S_{3}  \\
\beta S _{2}  & -S_{3} & 0
\end{array} \right) ,\quad
G^{\prime} =
\left ( \begin{array}{ccc}
0       & \omega_{3}  & -\omega_{2} \\
-\beta\omega_{3} & 0      & \omega_{1} \\
\beta\omega_{2}  & -\omega_{1} & 0
\end{array} \right) \eqno(75)
$$
Here ${\bf f}_{j}$ is the new moving trihedral of a curve or surface,
$$
{\bf f}_{1}^{2}=\beta = \pm 1, \quad {\bf f}_{2}^{2}={\bf f}^{2}_{3}=1 \eqno(76)
$$
and $\omega_{j}$ are already some functions of $S_{j}, S_{kx}, n$.
So we have the (1+1)-dimensional M-0 equation
$$
C^{\prime}_{t} - G^{\prime}_{x} + [C^{\prime}, G^{\prime}] = 0 \eqno(77)
$$

If remember that  $k, \tau,
 \sigma, \omega_{j}$ are some functions of $E, F, G, L, M, N $ (coefficients
of first and second fundamental forms) then the M-LXVI and M-0 equations (74)
and (77) describe some special class of surfaces with the some restriction
for $E, F, G, L, M, N$ which  follows from (71).
Let us rewrite the M-0 equation (77) in the form
$$
S_{t}-\frac{1}{n}V_{x}+[S,V]=0  \eqno(78)
$$
where
$$
V =
\frac{1}{2i}
\left ( \begin{array}{cc}
\omega_{1}             & \omega_{3}-i\omega_{2} \\
\omega_{3}+i\omega_{2} & - \omega_{1}
\end{array} \right) \eqno(79)
$$
The LR of the equation (78) is given by
$$
\psi_{x} = nS \psi \eqno(80a)
$$
$$
\psi_{t} =V\psi  \eqno(80b)
$$
If we take
$$
V=-in SS_{x} - 2in^{2}S \eqno(81)
$$
then from the equation (78) we obtain the isotropic
Landau-Lifshitz equation (LLE)
$$
2iS_{t} = [S, S_{xx}] \eqno(82)
$$
where $n$ plays the role of spectral parameter. Finally we consider the
gauge transformation
$$
\left ( \begin{array}{ccc}
{\bf f}_{1} \\
{\bf f}_{2} \\
{\bf f}_{3}
\end{array} \right)= E
\left ( \begin{array}{ccc}
{\bf e}_{1} \\
{\bf e}_{2} \\
{\bf e}_{3}
\end{array} \right)
\eqno(83)
$$
Then hence and from (1) and (74) we get
$$
C=E^{-1}C^{\prime}E - E^{-1}E_{x}, \quad
G=E^{-1}G^{\prime}E - E^{-1}E_{t}, \eqno(84)
$$
and the M-LXX equation give rise to the NLSE. In this case we speak that
the NLSE is L-equivalent (Lakshmanan equivalent) counterpart of the LLE [9].
As known between these equations take place G-equivalence (gauge
equivalence) [24].
The other reductions of the
M-LXVI and M-0 equations, including  and multi-dimensional spin systems and
for details,  see e.g. [22].

\section{Summary}

In this paper, we have considered some aspects of the relation
between differential geometry of curves and surfaces and multiidimensional
soliton equations. Our approach permits find some integrable classess
of curves and surfaces (soliton geometry). Starting from the known
fact that many soliton equations ,
for example, in 2+1 dimensions are particular cases of the mM-LXII
(and/or M-LXII) equation, we have shown they are exact reductions
of the SDYM equation. The connection between  spin systems and
curves/surfaces  is also discussed.  Finally we note that the more
detailed presentation of the above discussed  results were given
in our previous notes (see, e.g.  [18, 19, 22]).

\section{Acknowledgements}
RM would like to thank Prof. M. Lakshmanan for  hospitality
during visit,  for many stimulating discussions and for financial support.
RM also would like to thank Radha Balakrishnan and M. Daniel for hospitality
and for useful discussions.

\end{document}